%Paper: alg-geom/9602003
%From: Steven Bradlow <bradlow@vortex.math.uiuc.edu>
%Date: Fri, 2 Feb 96 13:58:23 -0600

%%%%%%%%%%%%%%%%%%%%%%%%%%%%%%%%%%%%%%%%%
%% A Hitchin-Kobayashi Correspondence %%%
%%        for Coherent Systems        %%%
%%         on Riemann Surfaces        %%%
%%                  by                %%%
%%%         Steven B. Bradlow         %%%
%%                 and                %%%
%%%         Oscar Garcia-Prada        %%%
%%%%%%%%%%%%%%%%%%%%%%%%%%%%%%%%%%%%%%%
%%%%%%%        AMSTeX           %%%%%%%
%%%%%%%%%%%%%%%%%%%%%%%%%%%%%%%%%%%%%%%

\input amstex
\documentstyle{amsppt}
\magnification=1200
\NoBlackBoxes
\TagsOnRight
\pageno=0
\refstyle {A}
\widestnumber\key{BDGW}

\def \qed{\hfill$\square $}
\def \-->{\longrightarrow}
\def \kler{K\"ahler}

\def \dbar{\overline{\partial}}
\def \dbare{\overline{\partial}_E}

\def \G{\frak G}
\def \Gc{\frak G_\Bbb C}

\def\CS{(\Cal E, V)}

\def\astable{$\alpha$-stable}
\def\adeg{\deg_{\alpha}}
\def\aslope{\mu_{\alpha}}
\def \bigspace{\Cal C\times(\Omega^0(E))^k}
\def \cs{\dbare,\phi_1,\dots,\phi_k}

\def \phik{\phi_1,\dots,\phi_k}

\def \GCcs{\Gc\times GL(k)}

\def \Gcs{\G\times \U(k)}

\def \X{\Cal X^k}
\def \Xcs{\Cal X^{CS}}

\def \E{\Cal E}
\def \C{\Cal C}
\def \I{\bold I}
\def\deg{{\operatorname{deg}}}
\def\rank{{\operatorname{rank}}}
\def\dim{{\operatorname{dim}}}
\def\Tr{{\operatorname{Tr}}}
\def\End{{\operatorname{End}}}

\def\Ker{{\operatorname{Ker}}}
\def\Im{{\operatorname{Im}}}
\def\max{{\operatorname{max}}}
\def\GL{{\operatorname{GL}}}
\def\U{{\operatorname{U}}}

%%%%%%%%%%%%%%%%%%%%%%%%%%%%%%%%%%%%%

\topmatter
\title
A Hitchin-Kobayashi Correspondence\\
for Coherent Systems\\
on Riemann Surfaces
\endtitle

\author Steven B. Bradlow \\
{\eightpoint Department of Mathematics \\
University of Illinois \\
Urbana, IL 61801 \\
bradlow\@uiuc.edu}
\\
\\
Oscar Garc\'{\i}a-Prada \\
{\eightpoint Departamento de Matematicas \\
Universidad Autonoma de Madrid \\
Ciudad Universitaria de Cantoblanco\\
 28049 Madrid,SPAIN\\
ogprada\@ccuam3.sdi.uam.es}
\\
\endauthor

\abstract
A `coherent system' $(\Cal E,V)$, consists of a holomorphic bundle plus a
linear subspace of its space of holomorphic sections. Based on the usual
notion in Geometric Invariant Theory, a notion of slope stability has been
defined for such objects (by Le Potier, and also by Rhagavendra and
Vishwanath). In this paper we show that stability in this sense is equivalent
to the existence of
solutions to a certain set of gauge theoretic equations. One of the equations
is essentially the vortex equation (i.e. the Hermitian-Einstein equation with
an additional zeroth order term), and the other is an orthonormality condition
on a frame for the subspace $V\subset H^0(\Cal E)$.

\endabstract
\endtopmatter
\leftheadtext{}
\rightheadtext{}
\bigskip
\eject
\document

%%%%%%%%%%%%%%%%%%%%%%%%%%%%%%%%%%%%%%%%%%%%%%%
\heading
\S 1. Introduction
\endheading
%%%%%%%%%%%%%%%%%%%%%%%%%%%%%%%%%%%%%%%%%%%%%%%

The Hitchin-Kobayashi correspondence for holomorphic bundles
relates
the existence of Hermitian-Einstein metrics to the property
of (slope)
stability.  Analagous correspondences are known to hold for
a number of
so-called augmented holomorphic bundles.  In this paper we
examine the
case of \it coherent systems\rm\ over Riemann surfaces.

Introduced in [LeP] and [RV], a coherent system is an `augmented bundle'
consisting of a holomorphic bundle together with a linear subspace of its
space
of holomorphic sections.  In order to construct moduli spaces for such
objects,
a notion of stability is required. As defined by LePotier, and also by
Rhagavendra and Vishwanath,
the definition involves a real parameter, which we label $\alpha$. Using this
notion, moduli spaces of $\alpha$-stable coherent systems are constructed in
[KN], [RV], and [LeP].  The constructions
given in these references are all based on geometric
invariant theory, and
relate the above notion of stability to standard notions of
stability in GIT. The parameter $\alpha$\ corresponds to a choice of
linearization for the
group action in the GIT setting.

In [BDGW] we introduced a set of geometric equations on
coherent
systems.  These equations, which we called the \it orthonormal vortex
equations\rm, are very similar to the vortex equations defined
in [B] on
holomorphic pairs and their generalization defined on $k$-
pairs in [BeDW].
For a given coherent system $(\Cal E,V)$\ the equations
determine a metric on $\E$\ and also a frame of $k$\ linearly independent
sections for $V\subset H^0(X,\Cal E)$.
In addition to the parameter $\alpha$, a second parameter (which we label
$\tau$) enters the equations. The two parameters are related by an identity
involving numerical invariants of $(\Cal E,V)$, viz. the rank and degree of
$\Cal E$\ and the dimension of $V$.

By a Hitchin-Kobayashi
corrrespondence for coherent systems, we mean a correspondence between the
property
$\alpha$-stability, and the existence of solutions to the orthonormal vortex
equations.  One direction in this correspondence was proven in [BDGW], namely

\proclaim{Theorem A [BDGW, Prop. 2.9]} Suppose that for some suitable choice
of
parameters $\alpha$\ and $\tau$, a coherent system $(\Cal E,V)$\ admits a
solution to the orthonormal vortex equations.
Then $(\Cal E,V)$\ is a direct sum of $\alpha$-stable coherent sustems.  If
$r$\ is coprime to either $d$\ or $k$, and if $\tau$\ is
generic, then
$(\Cal E,V)$\ is $\alpha$-stable
\endproclaim

In this paper we prove the converse, namely

\proclaim{Theorem B (3.13)} Fix $\alpha >0$, and let $(\Cal E,V)$\ be an
$\alpha$-stable coherent system. Then, with a suitably compatible choice for
$\tau$,
there is a unique smooth solution to the orthonormal vortex equations on
$(\Cal
E,V)$.
\endproclaim

The proof of Theorem A is a relatively minor modification of the proof in [B]
of the analogous statement for holomorphic pairs. That, in turn, is based on
the original arguement of Kobayashi showing that the Hermitian-Einstein
condition implies stability. The modifications required are described in
[BDGW].

The proof of Theorem B is the main result of this paper and is given in
Section 3. Before proceeding to this main result, we need to give some
definitions and
establish some basic  properties of $\alpha$-stability. This is done in
Section
2.

%%%%%%%%%%%%%%%%%%%%%%%%%%%%%%%%%%%%%%%%%%%%%%%
\heading
\S 2. Definitions and Basic properties
of Stable Coherent Systems
\endheading
%%%%%%%%%%%%%%%%%%%%%%%%%%%%%%%%%%%%%%%%%%%%%%%

\proclaim{Definition 2.1}Let $X$\ be a Riemann surface and
let
$E\-->X$\ be a fixed smooth complex bundle. A coherent
system on
$E$\ is a pair $(\Cal E,V)$\
where $\Cal E$\ is a holomorphic bundle isomorphic to $E$,
and $V$\ is a
linear subspace of $H^0(X,\Cal E)$. We say that $(\Cal
E,V)$\ is of type
(d,r,k) if $\deg(\Cal E)=d$, $\rank(\Cal E)=r$, and
$\dim(V)=k$.
\endproclaim

The equations we consider are equations for a metric on $E$\ and a basis
for $V$. We assume that $X$ has a fixed (\kler) metric.  Then, if we denote
a Hermitian bundle metric on $E$\  by $H$, and let
$\{\phi_1,\dots,\phi_k\}$\ be a set of linearly independent holomorphic
sections spanning $V$, the \it orthonormal vortex equations\rm\  can
be written as

$$\align
i\Lambda F_H+\Sigma_{i=1}^k\phi_i\otimes\phi_i^*&=\tau\bold
I\tag 2.1a\\
<\phi_i,\phi_j>&=\alpha\bold I_k\tag 2.1b
\endalign$$

In the first equation, $F_H$\ is the curvature of the metric
connection on
$\Cal E$, $\Lambda F_H$\ denotes the contraction with the
\kler\ form on
$X$, the adjoint in $\phi^*$\ is with respect to the
metric, $\bold
I$\ is the identity section of $\End E$, and $\tau$\ is a
real parameter.
The left hand side of the second equation is the $k\times
k$\ matrix
whose i-j entry is the $L^2$\ inner product  of
$\phi_i$\ and $\phi_j$\ in
$\Omega^0(X, E)$. On the right hand side, $\bold I_k$\
is the unit
$k\times k$\ matrix and $\alpha$\ is a real parameter.  By
taking
$\int_{X}\Tr$\ of the first equation, we see that
$\tau$\ and
$\alpha$\ must be related by
$$ d+\alpha k=r\tau\ .\tag 2.2$$

\proclaim{Definition 2.2}
\roster
\item Define the \bf
subobjects\it\  of
$\CS$ to
be subbundles $\Cal E'\subset\Cal E$\ together with
subspaces $V'\subset V\cap H^0(X,\Cal E')$.
\item A subpair is
\bf proper\it\ unless $\Cal E'=0$\ and $V'=0$, or
$\Cal E'=\Cal E$\ and $V'=V$.
\item A morphism $u: (\Cal E,V_E)\-->(\Cal F,W)$\ between
two coherent systems consists of a sheaf map $u:\Cal E\-->\Cal F$
such that $u(V_E)\subset V_F$, where here $u$\ denotes the induced map
on the space of global section.
\endroster
\endproclaim

\proclaim{Definition 2.3}
 For a given
$\alpha\in\Bbb R$, define the \bf $\alpha$-degree\it\  of
$(\Cal E',V')$\ to be
$$\adeg(\Cal E',V')=\deg(\Cal E')+\alpha \dim(V')\ .\tag 2.3$$
The \bf $\alpha$-slope\it\  of $(\Cal E',V')$\ is then
$$\aslope(\Cal E',V')=\frac{\adeg(\Cal E',V')}{\rank(\Cal
E')}\ .\tag 2.4$$
We say the Coherent Systems $\CS$\ is \bf
\astable\ \it if
for all subsystems $(\Cal E',V')$\
$$\aslope(\Cal E',V')<\aslope(\Cal E,V)\ .\tag 2.5$$
If the strong inequality is replaced by a weak one, then we
say the coherent system is $\alpha$-semistable.
\endproclaim

\proclaim{Definition 2.4} A coherent system $(\Cal E,V)$ is \bf
decomposable\it\ if the bundle decomposes as $\Cal E=\Cal
E_1\oplus\Cal
E_2$\ and $V=V_1\oplus V_2$, with $V_i\subset H^0(X,
\Cal E_i)$\
for $i=1,2$.  That is, a decomposable coherent system splits
as
$(\Cal E,V)=(\Cal E_1,V_1)\oplus (\Cal E_2,V_2)$.
\endproclaim

It follows immediately that
\proclaim{Lemma 2.5} Suppose that a type $(d,r,k)$\ coherent
system
decomposes into coherent systems of type $(d_1,r_1,k_1)$\
and
$(d_2,r_2,k_2)$, i.e.
$(\Cal E,V)=(\Cal E_1,V_1)\oplus (\Cal E_2,V_2)$.  Then for
any $\alpha$,
$$\mu_{\alpha}(\Cal E,V)=\frac{r_1}{r}\mu_{\alpha}(\Cal
E_1,V_1)+
\frac{r_2}{r}\mu_{\alpha}(\Cal E_2,V_2)\ .$$
\endproclaim

\proclaim{Definition 2.6} A coherent system $(\Cal E,V)$ \ is called \bf
polystable\it\ if it
is decomposable and decomposes into a sum of $\alpha$-stable
coherent
systems, each with $\alpha$-slope $\mu_{\alpha}(\Cal E,V)$.
\endproclaim

\proclaim{Proposition 2.7}
Let $(\Cal E,V)$ be a  coherent system of type $(d,r,k)$,  and let
$(\Cal E',V')$ be a subsystem of type $(d',r',k')$\ such that
$\mu_{\alpha}({\Cal E},V)=\mu_{\alpha}({\Cal E}',V')$. Then either
$$
rk'=r'k\ \text{and}\ \mu({\E}')=\mu(\E),
$$
or
$$
\alpha = \frac{r'd-rd'}{rk'-r'k}.
$$
In particular, if $r$\ is coprime to either $d$\ or $k$, and
$\alpha$ is not a rational number with denominator of
magnitude less
than or equal to $rk$, then all $\alpha$-semistable coherent systems are
$\alpha$-stable.
\endproclaim

For convenience, we will say that \it $\alpha$\ is
generic\rm\  if it is not
of the above form, i.e. if it is not a rational number with
denominator of
magnitude less than or equal to $rk$.  Similarly, values of
$\tau$\ which
correspond, via (1.5), to generic values of $\alpha$\ will be called generic.

It has been noted by ( by Lepotier and by King and Newstead) that under
the above definitions, coherent systems do not form a convenient category.
In particular, quotients of coherent systems are not necessarily included in
the category. (Consider, as an extreme example of this possibility, the
quotient of
$(\Cal E,V)$\ by $(\Cal E,W)$, where $W\subset V$.) Following King and
Newstead, one
thus needs to introduce a larger category, which we will denote by $\Cal
{CS}$.
 The objects in $\Cal {CS}$\ are triples $(\Cal E, \Bbb V, \sigma)$, where
$\Cal E$\ is any sheaf, $\Bbb V$\ is a finite dimensional vector space, and
$\rho: \Bbb V\otimes\Cal O_X\-->\Cal E$\ is a sheaf map (not necessarily
injective). A morphism between two such objects, say $(\Cal E, \Bbb V, \rho)$\
and
$(\Cal F, \Bbb W, \sigma)$, consists of  a sheaf map $f:\Cal E\-->\Cal F$\ and
a linear map $L:\Bbb V\-->\Bbb W$\
such that the folllowing diagram commutes
$$\CD
\Bbb V\otimes\Cal O_X@>\rho >> \Cal E\\
@V L\otimes 1VV  @VV fV\\
\Bbb W\otimes\Cal O_X@>\sigma >> \Cal F
\endCD$$
Then:
\roster
\item The coherent systems correspond to the objects $(\Cal E, \Bbb V, \rho)$\
in which
$\Cal E$\ is torsion free and $\rho$\ is injective. We will denote the
corresponding coherent system by $(\Cal E, V)$, where $V\subset H^0(X, \Cal
E)$\ is the image $\rho(\Bbb V\otimes\Cal O_X)$.
\item The notion of the type of a coherent system extends to objects in $\Cal
CS$.
\item Subobjects of $(\Cal E, \Bbb V, \rho)$\ are objects
$(\Cal E', \Bbb V', \rho')$\  with a morphism made up of inclusions $i_E:\Cal
E\-->\Cal E'$
and $i_V:\Bbb V\-->\Bbb V'$,
\item The definition of $\alpha$-stability can be extended
to objects in $\Cal CS$,  and one can show (cf [KN])
that the $\alpha$-semistable objects are in fact coherent systems.
\endroster

When considering semistable objects, we thus need not distinguish between
coherent systems and objects in the larger category $\Cal CS$. The larger
catgeory will however be needed at some places in the ensuing discussion,
especially in the construction of filtrations of coherent systems (cf. Lemma
2.9).

\proclaim{Lemma 2.8 }Fix $\alpha >0$, and let $(\Cal E_1,V_1)$\ and
$(\Cal E_2,V_2)$\ be $\alpha$-semistable coherent systems.
Let $u:(\Cal E_1,V_1)\-->(\Cal E_2,V_2)$\ be a morphism of
coherent systems.  Suppose that $\mu_{\alpha}(\Cal
E_1,V_1)\ge\mu_{\alpha}(\Cal
E_2,V_2)$\ and that in the case of equality, at least one of
the coherent systems is $\alpha$-stable.

\roster
\item If the induced map $u:V_1\--> V_2$\ is injective, then
 $u:\Cal E_1\-->\Cal E_2$\ is an isomorphism.
%either the intersection $V_1\cap H^0(\Ker (u))$\ is
%non-empty, or $u$\ is an isomorphism,

\item If the induced map $u:V_1\--> V_2$\ is the zero map, then
$u:\Cal E_1\-->\Cal E_2$\ is the zero map.

%if $V_1\subset H^0(\Ker(u))$, then $u=0$.

\endroster
\endproclaim

\demo{Proof}

(1) Suppose that $u:V_1\-->H^0(\Cal E_2)$\ is injective. Consider the
subsytems
$(\Ker(u), 0)$\ and $(\Im(u),u(V_1))$\ of $(\Cal E_1,V_1)$\
and $(\Cal E_2,V_2)$\ respectively.  By the injectivity assumption, we have

$$0\-->(\Ker(u), 0)\-->(\Cal E_1,V_1)\-->(\Im(u), u(V_1))
\-->0\ .$$

If $\Ker(u)\ne O$, then the $\alpha$-semistability condition implies
$$\mu_{\alpha}(\Ker(u), 0)\le \mu_{\alpha}(\Cal E_1,V_1)\ ,
\tag 2.6$$
$$\mu_{\alpha}(\Im(u), u(V_1))\le \mu_{\alpha}(\Cal E_2,V_2)\
.\tag 2.7$$
Moreover, the inequality (2.6) (respectively (2.7)) is strict if $(\Cal
E_1,V_1)$, (respectively $(\Cal E_2,V_2)$) is $\alpha$-stable.  Also,
$$r_k\mu_{\alpha}(\Ker(u), 0)+r_I\mu_{\alpha}(\Im(u), u(V_1))=
r_1\mu_{\alpha}(\Cal E_1,V_1)\ ,\tag 2.8$$
where $r_k$\ and $r_I$\ denote the ranks of $\Ker(u)$\ and
$\Im(u)$\ respectively.  We thus find that either
$$\mu_{\alpha}(\Cal E_1,V_1)>\mu_{\alpha}(\Cal
E_2,V_2)\ge\mu_{\alpha}(\Cal E_1,V_1)\ ,$$
or
$$\mu_{\alpha}(\Cal E_1,V_1)=\mu_{\alpha}(\Cal
E_2,V_2)>\mu_{\alpha}(\Cal E_1,V_1)\ .$$
Thus $\Ker(u)=O$\ and $u$\ is an isomorphism.

(2) If $V_1\subset H^0(\Ker (u))$, then we can consider the
subsystems $(\Ker(u), V_1)$\ and $(\Im(u),0)$. Part (2) then
follows by similar arguements to the ones above.

\hfill$\square $
\enddemo

\noindent\bf Remark:\rm\ It follows from part (1) that,
under the assumptions of Lemma 2.8,
$u:\Cal E_1\-->\Cal E_2$\ is an isomorphism if and only if
the induced map $u:V_1\--> V_2$\ is an isomorphism.

\proclaim{Lemma 2.9} (Harder Narasimhan filtration for
unstable coherent systems)
Fix $\alpha >0$, and let $(\Cal E,\Bbb V,\rho)$\ be an object in $\Cal {CS}$
which is not $\alpha$-semistable.  Then there is a unique
filtration by sub-coherent systems

$$0=(\Cal E_0,V_0)\subset(\Cal E_1,V_1)\subset\dots(\Cal
E_k,V_k) \ ,\tag 2.9$$
such that
\roster
\item $\Cal E_k=\Cal E$,
\item if $(\Cal E,\Bbb V,\rho)$\ corresponds to a coherent system (i.e. if
$\rho$\ is injective)
then $V_k=V$. Otherwise, $V_k=\rho(\Bbb V\otimes \Cal O_X)$
\item for $1\le i\le k$, $(\Cal E_i/\Cal E_{i-1}, V_i/V_{i-1})$\ is an
$\alpha$-semistable
sub-system of $(\Cal E/\Cal E_{i-1},V/V_{i-1})$,
\item  the $\alpha$-slopes are ordered
such that
$$\mu_{\alpha}(\Cal E_1,V_1)>\mu_{\alpha}(\Cal E_2/\Cal E_1,
V_2/V_1)>\dots\mu_{\alpha}(\Cal
E_k/\Cal E_{k-1}, V_k/V_{k-1})\ .\tag 2.10$$
\endroster
\endproclaim

\demo{Proof} As for holomorphic bundles, the result follows
from the fact that the slopes of subobjects are bounded
above. If $\mu(\Cal E')\le\mu_{\max}$\ for all
subbundles, then for any subobject $(\Cal E',\Bbb V',\rho')$\ we have
$$\frac{\deg(\Cal E') +\alpha \dim(\Bbb V')}{\rank(\Cal E')}\le
\mu_{\max}+\alpha
\frac{k}{r}\
,$$
where $k=\dim(\Bbb V)$\ and $r=\rank(E)$. It follows that there are
subobjects of maximal $\alpha$-slope, and that amongst these there is one, say
$(\Cal E_1,\Bbb V_1,\rho_1)$, of maximal rank. By construction, this maximal
subobject is
$\alpha$-semistable, and thus corresponds to a sub-coherent system, which we
denote by $(\Cal E_1,V_1)$.
Also, by the maximality property of this sub-object, we can assume that
$V_1=V\cap H^0(\Cal E_1)$.

The rest of the filtration is constructed by iterating this construction.

\qed
\enddemo

\proclaim{Lemma 2.10}[RV, Cor 1.8] (Seshadri filtration for $\alpha$-
semistable coherent systems)

Fix $\alpha$, and let $(\Cal E,V)$\ be a coherent system
which is $\alpha$-semistable.  Then there is a filtration by sub-coherent
systems

$$0=(\Cal E_0,V_0)\subset(\Cal E_1,V_1)\subset\dots(\Cal
E_k,V_k)=(\Cal E,V)\tag 2.11$$
such that
\roster
\item $\mu_{\alpha}(\Cal E_i/\Cal E_{i-1}, V_i/V_{i-1})
=\mu_{\alpha}(\Cal E,V)$, and
\item $(\Cal E_i/\Cal E_{i-1}, V_i/V_{i-1})$\ is an $\alpha$-stable pair for
all $i$.
\endroster
Moreover, the isomorphism class of the pair
$$gr_{\alpha}(\Cal E,V)=\bigoplus(\Cal E_i/\Cal E_{i-1}, V_i/V_{i-1})$$
is independent of the filtration.
\endproclaim

%%%%%%%%%%%%%%%%%%%%%%%%%%%%%%%%%%%%%%%%%%%%%%%
\heading
\S 3. The main result
\endheading
%%%%%%%%%%%%%%%%%%%%%%%%%%%%%%%%%%%%%%%%%%%%%%%

In this section we prove Theorem B.  The proof follows the method introduced
by
Donaldson in [Do], and later used in [H] and in [GP]. Thus we reformulate the
vortex equations as minimization criteria for a functional based on a
symplectic moment map.

We begin with some background and notation. Let  $E\-->X$\ be the underlying
smooth bundle for $\E$, i.e. $\E$\ is $E$\ together with a holomorphic
structure.
Specification of a holomorphic structure is equivalent to a
choice of a $\dbar$-operator
on $\Omega^0(E)$. We denote the set of all such operators by
$\C$.
A dimension-$k$ coherent system on $E$\ consists of a pair
$(\E,V)$\ where $V$\ is a $k$-dimensional linear subspace of
$H^0(X,\E)$.   If we use $k$-frames in $H^0(X,\E)$\ to
describe such subspaces, then the configuration space of all
such coherent systems can be described as follows.
Define:

$$\X=\bigspace\ ,\tag 3.1$$

$$\X_{0} = \{(\dbare,\phik)\in \X : \text
{the sections are linearly independent}\}\ ,\tag 3.2$$

$$\Cal H^k = \{(\dbare,\phik)\in \X :
\dbare(\phi_i)=0\ \text
{for }\ i=1,\dots k \}\ ,\tag 3.3$$

$$ \Xcs = \X_{0}/\GL(k), \text{where}\  \GL(k)\ \text{acts
on}\  (\Omega^0(E))^k\ .\tag 3.4$$

\noindent The space of  dimension-$k$
coherent systems is then the subvariety of $\Xcs$\ given by

$$\Cal H^{CS}=(\X_{0}\cap\Cal H^k)/\GL(k)\ .\tag 3.5$$

Suppose now that we fix a Hermitian metric on $E$.  Then we can use this fixed
bundle metric to define the
unitary gauge group $\G$.  Also the spaces $\Cal C$\ and
$\Omega^0(E)$\ acquire \kler, and thus symplectic,
structures.
Denoting the symplectic forms by $\omega_{\Cal C}$\ and
$\omega_0$\ respectively, we can give  $\bigspace$\ the
symplectic structure with symplectic form
$\omega=\omega_{\Cal C}
+\omega_0+\dots +\omega_0$.
The gauge group $\G$\ acts by

$$g(\dbare,\phik)=(g\circ\dbare\circ
g^{-1},g\phi_1,\dots,g\phi_k)\ \tag 3.6$$

and this action preserves $\omega $.  In addition, the
action of
$\U(k)$\ on $(\Omega^0(E))^k$ is also symplectic, and
commutes
with the action of $\G$.

\proclaim{Proposition 3.1}There are moment maps for the actions
of $\G$\ and $\U(k)$
on $\X$.  These are given, respectively,  by
$$\align
 \Psi_{\G}(\dbare,\phik)&=\Lambda F,\\
 \Psi_\U(\dbare,\phik)&=-i<\phi_i,\phi_j>
\endalign$$
Here $F$\ denotes the curvature of the metric connection determined by
$\dbare$\ and the fixed Hermitian metric on $E$.
\endproclaim

Notice that for any non-zero real number $\alpha$, the union
of the
$\GL(k)$-orbits through the level set
$\Psi_U^{-1}(-i\alpha\bold I_k)$\ is exactly $\X_{0}$. It
follows from
this that $\Xcs$\ can be described as the symplectic
reduction
$\Psi_U^{-1}(-i\alpha\bold I_k)/\U(k)$. Similarly,
$\Cal H^{CS}= (\Psi_U^{-1}(-i\alpha\bold I_k)\cap \Cal
H^K)/\U(k)$.

The actions (on $\X$) of $\Gc$\ and $\GL(k)$\ commute,  and
thus the $\Gc$\ action descends
to the quotient $\Xcs = \X_{0}/\GL(k)$. The subvariety $\Cal
H^{CS}$\ is preserved by this action. Similarly, in the
symplectic descriptions, $\G$
acts on $\Psi_U^{-1}(-i\alpha\bold I_k)/\U(k)$\ and this
action is symplectic with respect
to the reduced symplectic structure.  Furthermore, there is
a moment map for this action
given as follows.

$$\Psi^{CS}_{\G}(\Cal E,V)=\Lambda F-
i\Sigma\phi_i\otimes\phi_i^*\tag 3.7$$
where $\{\phi_i,\dots,\phi_k\}$\ is any frame for $V$\ such
that $<\phi_i,\phi_j>=\alpha\bold I_k$.

Theorem B can now be reformulated as

\proclaim{Theorem 3.2} If $(\Cal E,V)$\ is an $\alpha$-stable
coherent
system, then
there is a unique smooth solution to
$$\Psi^{CS}_{\G}(\Cal E,V)=-i\tau\bold I$$
on the $\Gc$\ orbit through $(\Cal E,V)$.
\endproclaim

Following [Do], we use the moment map (3.7) to define an analog of the
functional
$J$\ defined in [Do]. As in
[Do] we use $\nu$\ to denote the Trace norm on matrices, and
set
$$N(A)=(\int \nu(A)^2)^{1/2}\ .\tag 3.8$$

\proclaim{Definition 3.3}With $\tau \rank(E)=\deg(E)+\alpha k$, define
$J:\Xcs\-->\Bbb R$\ by
$$J(\Cal E,V)=N(\Psi^{CS}_{\G}(\Cal E,V)+i\tau\I)\ .\tag 3.9$$
\endproclaim

The basic idea of the proof of Theorem 3.2 is to study the restriction of $J$\
to the $\Gc$-orbit through
$(\E,V)$. The main steps involved are
\roster
\item to show that a minimizing sequence on this orbit converges, say to
$(\Cal
E_{\infty},V_{\infty})$,
\item to show that $(\Cal E_{\infty},V_{\infty})$\ is on the orbit through
$(\Cal E,V)$, and
\item to show that $J=0$\ at this minimizer, i.e. $J(\Cal
E_{\infty},V_{\infty})=0$.
\endroster

\bigskip

%%%%%%%%%%%%%%%%%%%%%%%%%%%%
\subheading{\hfil Step 1: Convergence of a minimizing sequence}
%%%%%%%%%%%%%%%%%%%%%%%%%%%%

Let $\Cal O(\E,V)$\ be the $\Gc$-orbit through $(\E,V)$, and let
$\{(\E_n,V_n)\}$\ be a minimizing sequence for $J|_{\Cal
O(\E,V)}$. Let $(\E_n,V_n)$\
be represented by $(\dbar_n, \phi^n_1,\dots,\phi^n_k)$\ in
$\Psi_U^{-1}(-i\alpha\bold I_k)$. Then for each $n$\ we have

$$J(\E_n,V_n)=N(\Lambda F_n -
i\Sigma\phi^n_i\otimes{\phi^n_i}^*+i\tau\I)\ ,\tag 3.10a$$
$$<\phi^n_i,\phi^n_j>=\alpha\bold I_k\ ,\tag 3.10b$$

\noindent where the notation $F_n$\ refers to the curvature of metric
connection determined by the fixed Hermitian metric on $E$\ and the
holomorphic
structure $\dbar_n$.

\proclaim{Proposition 3.4}
There exist $(\dbar_{\infty},\phi^{\infty}_1,\dots,
\phi^{\infty}_k)$\ in $L^2_1(\bigspace)$\ such that after
passing to a subsequence, and
up to equivalence under $\Gcs$,
\roster
\item $\dbar_n\rightharpoonup\dbar_{\infty}$\ in $L^2_1$,
\item $\phi^n_i\rightharpoonup\phi^{\infty}_i$\ in $L^2_1$.
\endroster
(Here we may assume that the unitary
connection corresponding to $\dbar_0$\ is used to define the
$L^2_1$\
norm.)
Furthermore, $\phi^{\infty}_i\ne 0$\ and
$\dbar_{\infty}\phi^{\infty}_i=0$\ for $i=1,\dots,k$.  Also,
$<\phi^{\infty}_i,\phi^{\infty}_j>=\alpha\bold I_k$.
\endproclaim

\demo{Proof}
The norm used in the definition of $J$\ is
equivalent
to the usual $L^2$\ norm.  We may thus assume that there is
a uniform
bound on $||\Lambda F_n -
i\Sigma\phi^n_i\otimes{\phi^n_i}^*+i\tau\I||^2_{L^2}$.  But
this is
equivalent to a uniform bound on
$$\Cal {YMH}_{\tau}((\dbar_n,\phi^n_1,\dots,\phi^n_k)=
||F_n||^2_{L^2}+\Sigma||D_n\phi^n_i||^2_{L^2} +
||\Sigma\phi^n_i\otimes{\phi^n_i}^*-\tau\I||^2_{L^2}\ ,$$
where $F_n=D_n^2$, i.e. $D_n$\ is the metric connection determined by the
fixed
Hermitian metric on $E$\ and the holomorphic structure $\dbar_n$.
Indeed, since $\dbar_n\phi^n_i=0$\ for all $i$, the two
functionals
differ only by a constant determined by the topology of the
bundle $E$
(cf. [B]).  Thus there are uniform bounds on each
of
$||F_n||^2_{L^2}$,    $||D_n\phi^n_i||^2_{L^2}$\ and
$||\Sigma\phi^n_i\otimes{\phi^n_i}^*-\tau\I||^2_{L^2}$.

 Since the base manifold $X$\ is a Riemann surface, it
follows by a
theorem of Uhlenbeck (cf. [U])
that the bound on $||F_n||^2_{L^2}$\ is enough to ensure the
weak convergence in $L^2_1$\ norm  of the $D_n$. Since $\dbar_n$\ is the
antiholomorphic part of $D_n$, this in turn leads to the weak convergence of
the $\dbar_n$.

Since $||\phi^n_i||^2_{L^2}=\alpha$\ for all $n$\ and $i$,
we may assume
that for fixed $i$, $\{\phi^n_i\}$\ converges weakly in
$L^2$. Furthermore,
by using the bound on  $||D_n\phi^n_i||^2_{L^2}$\ and the
convergence of
$\dbar_n$, we obtain a bound on $||D_0\phi^n_i||^2_{L^2}$.
The weak
convergence of $\{\phi^n_i\}$\ may thus be taken to be in
$L^2_1$.
We will denote the limit point by
$\{(\dbar_{\infty},\phi^{\infty}_1,\dots,
\phi^{\infty}_k)\}$.

By construction, we have $\dbar_n\phi^n_i=0$\ for each $n$\
and $i$.  Thus
$$\dbar_{\infty}\phi^{\infty}_i=(\dbar_{\infty}-
\dbar_n)\phi^{\infty}_i+
(\dbar_n-\dbar_0)(\phi^{\infty}_i-
\phi^n_i)+\dbar_0(\phi^{\infty}_i-\phi^n_i)\ .$$
Defining $x_n,\ y_n\in\Omega^{0,1}(\End(E))$\ by
$\dbar_{\infty}-\dbar_n=x_n$\ and $\dbar_n-\dbar_0=y_n$, we
get
$$||\dbar_{\infty}\phi^{\infty}_i||_{L^2}\le
||x_n||_{L^4}||\phi^{\infty}_i||_{L^4}+
||y_n||_{L^4}||\phi^{\infty}_i-\phi^n_i||_{L^4}+
||\phi^{\infty}_i-\phi^n_i||_{L^2_1}\ .$$
But, since $L^2_1\subset L^4_0$\ is compact, we can assume
$||x_n||_{L^4}\-->0$,
$||\phi^{\infty}_i-\phi^n_i||_{L^4}\-->0$, and also that
$||y_n||_{L^4}$\ is bounded. Thus
$\dbar_{\infty}\phi^{\infty}_i=0$.

Examining $<\phi_i,\phi_j>$, we find
$$\align
|<\phi^{\infty}_i,\phi^{\infty}_j>-\alpha\bold I_k|&\le
|<\phi^{\infty}_i,\phi^{\infty}_j>-<\phi^n_i,\phi^n_j>|+
|<\phi^n_i,\phi^n_j>-\alpha\bold I_k|\\
&\le|<\phi^{\infty}_i-\phi^n_i,\phi^{\infty}_j>-
<\phi^n_i,\phi^{\infty}_j-\phi^n_j>|\\
&\le\sum_{i,j}{||\phi^{\infty}_i-
\phi^n_i||_{L^2}||\phi^{\infty}_j||_{L^2}}+
\sum_{i,j}{||\phi^{\infty}_j-
\phi^n_j||_{L^2}||\phi^n_j||_{L^2}}\\
&\le k\alpha\sum_i{||\phi^{\infty}_j-\phi^n_j||_{L^2}}
\endalign$$
Thus $<\phi^{\infty}_i,\phi^{\infty}_j>=\alpha\bold I_k$. In
particular,
$\phi^{\infty}_i\ne 0$\ for all $i$.

\qed\enddemo

\proclaim{Definition 3.5} Let $\E_{\infty}$\ be the holomorphic
bundle determined
by $\dbar_{\infty}$, and let $V_{\infty}$\ be the subspace
of $H^0(X,\E_{\infty})$\
spanned by $\{\phi^{\infty}_1,\dots, \phi^{\infty}_k\}$.
Notice that because of (c), the dimension of $V_{\infty}$\
is k.
Thus $\{(\dbar_{\infty},\phi^{\infty}_1,\dots,
\phi^{\infty}_k)\}$\ defines a coherent system
$(\E_{\infty}, V_{\infty})$\ in $\Cal H^{CS}$.
\endproclaim

\bigskip

%%%%%%%%%%%%%%%%%%%%%%%%%%%%
\subheading{\hfil Step 2: the minimizer is on the orbit}
%%%%%%%%%%%%%%%%%%%%%%%%%%%%

We now show that if $(\E, V)$\ is $\alpha$-stable,
then $(\E_{\infty}, V_{\infty})$\ is on the same $\Gc$-orbit as $(\E,V)$.  We
do this by first showing that there is a non-trivial homomorphism
$h:(\E,V)\-->(\E_{\infty}, V_{\infty})$, and then proving
that
$h$\ must be an isomorphism if $(\E,V)$\ is $\alpha$-stable.

Since all  the $(\dbar_n, \phi^n_1,\dots,\phi^n_k)$\ in the minimizing
sequence
lie on the same $\GCcs$-orbit
in $\X$, for each $n$\ we can find $g_n\in\Gc$\ and
$A^{(n)}\in \GL(k)$\ such that
$\dbar_n=g_n(\dbar_0)$\ and
$A^{(n)}_{ij}\phi^n_j=g_n(\phi^0_i)$.  Let $\tilde{g}_n=g_n/||g_n||_{L^2}$,
and
set
$\tilde{A}^{(n)}=A^{(n)}/|A^{(n)}|$.  Since
$\{\tilde{A}^{(n)}\}$\ is a
bounded sequence of $k\times k$\ matrices, it has a
convergent subsequence. Say
$$\tilde{A}^{(n)}\-->\tilde{A}^{\infty}\ .$$

\proclaim{Lemma 3.6} The sequence $\{\tilde{g}_n\}$\ has a
subsequence which converges in $L^2_1$, say
$$\tilde{g}_n\-->\tilde{g}_{\infty} $$
Furthermore, $\tilde{g}_{\infty} \ne 0$,
$\dbar_{\infty}\circ\tilde{g}_{\infty}=\tilde{g}_{\infty}
\circ\dbar_0$,
and $\tilde{g}_{\infty}(V)\subset V_{\infty}$.  In
particular, $\tilde{g}_{\infty}$\
gives a nontrivial homomorphism from
$(\E,V)$\ to $(\E_{\infty}, V_{\infty})$.
\endproclaim

\demo{Proof}
The following arguement showing the convergence of
$\{\tilde{g}_ n\}$\ can be found in [H]. Using the
operators $\dbar_n$\ and $\dbar_0$\ we define an elliptic
operator
$$\dbar_{n,0}:\Omega^0(X, E\otimes E^*)\-->\Omega^{0,1}(X,
E\otimes E^*)\ .$$
We can write $\dbar_{n,0}=\dbar_{\infty,0}+\beta_n$, with
$\beta_n\rightharpoonup 0$
in $L^2_1$. From the fact that
$$\dbar_{n,0}(\tilde{g}_n)=\dbar_{n,0}(g_n)=0\ ,$$
we get
$$||\dbar_{\infty,0}(\tilde{g}_n)||_{L^2}\le
||\beta_n||_{L^4}||\tilde{g}_n||_{L^4}\ .\tag 3.11$$
Thus,  elliptic estimates give
$$||\tilde{g}_n||_{L^2_1}\le
C(||\beta_n||_{L^4}||\tilde{g}_n||_{L^4}+
||\tilde{g}_n||_{L^2})\ .\tag 3.12$$
But $||\tilde{g}_n||_{L^2}=1$. Also, since since
$L^2_1\subset L^4_0$\
is compact, we have $||\beta_n||_{L^4}\-->0$\ and thus (3.12)
gives a uniform bound on $||\tilde{g}_n||_{L^2_1}$.  After
renaming a subsequence, we thus get
$\tilde{g}_n\rightharpoonup\tilde{g}_{\infty}$\ in $L^2_1$.
Since $||\tilde{g}_n||_{L^2}=1$, we can conclude that
$\tilde{g}_{\infty}\ne 0$. It now
follows from (3.11) that
$\dbar_{\infty,0}(\tilde{g}_{\infty})=0$,
i.e.
$\dbar_{\infty}\circ\tilde{g}_{\infty}=\tilde{g}_{\infty}
\circ\dbar_0$.

 Finally, if we set $c_n=|A^{(n)}|/||g_n||_{L^2}$, we can
write
$$\tilde{g}_n(\phi^0_i)=c_n\tilde{A}^{(n)}_{ij}\phi^n_j\ .$$
Thus for each $n$, $\tilde{g}_n(\phi^0_i)$\ is in the space
spanned by $\{\phi^n_j\}^k_{j=1}$. If we define
$$d_n=\sum_{p=1}^{k}{\frac{|<\tilde{g}_n(\phi^0_i),\phi^n_p>
|^2}{||\phi^n_p||^2}
-||\tilde{g}_n(\phi^0_i)||^2}\ ,$$
so $d_n$\ is the distance from $\tilde{g}_n(\phi^0_i)$\ to
the space spanned
by $\{\phi^n_j\}^k_{j=1}$, then clearly $d_n=0$\ for all
$n$. Thus, using the convergence
of the $\phi^n_i$\ and the $\tilde{g}_n$, we get
$$0=\lim_{n\rightarrow \infty}{d_n}
=\sum_{p=1}^{k}{\frac{|<\tilde{g}_{\infty}(\phi^0_i),\phi^{\
infty}_p>|^2}
{||\phi^{\infty}_p||^2}-||\tilde{g}_{\infty}(\phi^0_i)||^2}\
.$$
That is, $\tilde{g}_{\infty}(\phi^0_i)$\ is in the space
spanned by $\{\phi^{\infty}_j\}^k_{j=1}$.

\qed\enddemo

It remains to establish that any such homomorphism is an
isomorphism.  As in [Do], the proof in the general case will involve induction
on the rank of the bundle.  We first observe that the result is true for rank
one bundles:

\proclaim{Lemma 3.7} If the rank of $E$\ is one, then
$\tilde{g}_{\infty}$\ is a
constant multiple of the identity. In particular,
$\tilde{g}_{\infty}:(\Cal E,V)\-->({\Cal
E}_{\infty},V_{\infty})$ is an
isomorphism.
\endproclaim

\demo{Proof} Since $\Cal E$\ and ${\Cal E}_{\infty}$\ are
line bundles of
equal
degree, any non trivial homomorphism
$\tilde{g}_{\infty}:\Cal E\-->{\Cal
E}_{\infty}$\ is a constant multiple of the identity.
\qed\enddemo

Suppose now that $\rank(E)>1$\ but that $\tilde g_{\infty}:(\Cal E,V)\-->(\Cal
E_{\infty},V_{\infty})$\ is not an isomorphism. Since
$\tilde g_{\infty}(V)\subset V_{\infty}$, the canonical
factorization of $\tilde g_{\infty}:\Cal E\-->\Cal
E_{\infty}$\ (cf. [Do]) gives a factorization
$$\CD
0@>>>(\Cal K,V_K)@>>>(\Cal E,V)@>>>(\Cal I,V/V_k)@>>>0\\
@.@.@V\tilde g_{\infty}VV@VVV\\
0@<<<(\Cal Q,V_Q)@<<<(\Cal
E_{\infty},V_{\infty})@<<<(\tilde{\Cal I},\tilde
g_{\infty}(V))@<<<0
\endCD\tag 3.13$$

In this diagram, $\rank(\Cal I)=\rank(\tilde{\Cal I})$,
$\dim(V/V_k)=\dim(\tilde g_{\infty}(V))$, $\deg(\Cal I)\le
\deg(\tilde{\Cal I})$, while $\rank(\Cal E)=\rank(\Cal
E_{\infty})$, $\dim(V)=\dim(V_{\infty})$, and $\deg(\Cal E)=
\deg(\Cal E_{\infty})$. Furthermore, since $(\Cal E,V)$\ is
\astable, we have $\mu_{\alpha}(\Cal
K,V_K)<\mu_{\alpha}(\Cal E,V)<\mu_{\alpha}(\Cal I,V/V_k)$.

Strictly speaking, the quotients $(\Cal I,V/V_k)$\ and $(\Cal Q,V_Q)$\ may
not
be coherent systems, and the above factorization makes sense only in the
category $\Cal CS$. In fact, $(\Cal I,V/V_k)$\ \it is\rm\ a coherent system,
i.e. the sheaf map in the corresponding object in $\Cal CS$\ is indeed
injective. (This is because $V_K=V\cap H^0(X,\Cal K)$. )  The only abuse of
notation is thus in the term $(\Cal Q,V_Q)$. We can tolerate this since we
will
not make use of this term in what follows.

We now need the following Lemmas, which are the analogs of Lemmas 2 and 3 in
[Do].

\proclaim{Lemma 3.8}Suppose that $\rank(\Cal E)>1$, and that $(\Cal E,V)$\ has
a
subsystem $(\Cal E_1,V_1)$\ with
$\mu_{\alpha}(\Cal E_1,V_1)\ge\mu_{\alpha}(\Cal E,V)$.
Set $\Cal E_2=\Cal E/\Cal E_1$\ and $V_2=V/V_1$, and let
$$\mu_{\alpha}(\Cal E_2,V_2)=
\mu(\Cal E_2)+\alpha\frac{\dim(V_2)}{\rank(\Cal E_2)}\ .$$

Let $\{\phik\}$\ be any basis for $V$\ such that the first
$k_1$\ sections form a basis for $V_1$.

Then
$$J(\cs)\ge
 R_1(\mu_{\alpha}(\Cal E_1,V_1)-\mu_{\alpha}(\Cal E,V))
+R_2(\mu_{\alpha}(\Cal E ,V) -\mu_{\alpha}(\Cal E_2, V_2))
\ .\tag 3.14$$
\endproclaim

\demo{Proof}

Let $(\Cal E_1,V_1)$\ be a subsytem of $(\Cal E,V)$. Let
$\{\phi_1,\dots,\phi_k\}$\ be an orthonormal frame for $V$\
such that $V'$\ is spanned by the first $k_1$\ elements,
$\{\phi_1,\dots,\phi_{k_1}\}$. With respect to the smooth
splitting $E=E_1\oplus(E/E_1)$, we get a decomposition of
each section into its components in $E_1$\ and
$E_2=(E/E_1)$, i.e.
$$\phi_i=\phi'_i + \phi^{\perp}_i\ .$$
The expression $\Sigma\phi_i\otimes\phi_i^*$\ thus has a
block decomposition as
$$\Sigma\phi_i\otimes\phi_i^*=\pmatrix
\sum_{i=1}^{k_1}{\phi_i\otimes\phi_i^*}+
\sum_{i=k_1+1}^{k }{\phi'_i\otimes(\phi'_i)^*}& * \\
* & \sum_{i=k_1+1}^{k
}{\phi^{\perp}_i\otimes(\phi^{\perp}_i)^*}
\endpmatrix $$

\noindent Notice that we can write
$$\int\{\Tr(\sum_{i=1}^{k_1}{\phi_i\otimes\phi_i^*}+
\sum_{i=k_1+1}^{k }{\phi'_i\otimes(\phi'_i)^*})\}=k_1 +
\int{t^2}\ ,$$
where $t^2$\ is a non-negative real number. It follows that
$$\int\{\Tr(\sum_{i=k_1+1}^{k
}{\phi^{\perp}_i\otimes(\phi^{\perp}_i)^*})\}=
(k-k_1)-\int{t^2}\ .$$

\noindent Recall also the block decomposition of $i\Lambda
F_H$\ , viz.
$$i\Lambda F_H=\pmatrix
i\Lambda F_{H_1}+\Pi& * \\
* & i\Lambda F_{H_2}-\Pi
\endpmatrix \ ,$$
where $\Pi$\ is a positive definite endomorphism coming from
the second fundamental form of the inclusion of $E_1$\ in
$E$.

By exactly the same arguement as in Lemma 2 of [Do], we thus get

$$\align \nu(i\Lambda F_H+\Sigma\phi\otimes\phi^*-
\tau\bold I)&\ge
|\Tr(i\Lambda F_{H_1}+
\sum_{i=1}^{k_1}{\phi_i\otimes\phi_i^*-
\tau\bold I_1}+\Tr\Pi+t^2|\\
&+|\Tr(i\Lambda F_{H_2}+
\sum_{i=K_1+1}^{k}{\phi^{\perp}_i\otimes
(\phi_i^{\perp})^*-
\tau\bold I_1}-\Tr\Pi|\ ,
\endalign$$
But $\tau=\mu_{\alpha}(\Cal E,V)$\ and by assumption
$\mu_{\alpha}(\Cal E_1,V_1)\ge\mu_{\alpha}(\Cal E,V)\ge
\mu_{\alpha}(\Cal E_2,V_2) $. Thus we get
$$J(\cs) \ge
 R_1(\mu_{\alpha}(\Cal E_1,V_1)-\tau) +
R_2(\mu_{\alpha}(\Cal E ,V) -\tau)\ .
$$
\qed\enddemo

\proclaim{Lemma 3.9}Let $(\Cal E,V)$\ be an $\alpha$-stable
coherent system, with $\rank(\Cal E)>1$.  Suppose that $(\Cal E,V)$\ is given
as
an
extension of coherent systems
$$0\-->(\Cal E_1,V_1)\-->(\Cal E,V)\-->(\Cal E_2,V_2)
\-->0\ $$
and suppose that Theorem 3.2 is true for coherent systems on bundles of lower
rank than $\Cal E$.
Let $\{\phik\}$\ be any basis for $V$\ such that the first
$k_1$\ sections form a basis for $V_1$.
Then
$$J(\dbare',\phi'_1,\dots,\phi'_k)\le
 R_1(\mu_{\alpha}(\Cal E ,V )-\mu_{\alpha}(\Cal E_1, V_1))
+R_2(\mu_{\alpha}(\Cal E_2 ,V_2) -\mu_{\alpha}(\Cal E, V))
\ ,\tag 3.15$$
for some $(\dbare',\phi'_1,\dots,\phi'_k)$\ on the
$\Gc\times \GL(k,\Bbb C)$\ orbit through $(\cs)$.
\endproclaim

\demo{Proof} Suppose first that $(\Cal E_1,V_1)$\ and $(\Cal E_2,V_2)$\ are
both $\alpha$-stable.  Let $\dbar_i$\ denote
the holomorphic structures on $\Cal E_i$. By our inductive
hypothesis, and possibly after a complex gauge transformation of $\Cal E_1$,
we
can
pick $\dbar_1$\ and
$\phi_1,\dots,\phi_{k_1}$\ such that
$$\align
&i\Lambda F_1+\Sigma_{i=1}^{k_1}\phi_i\otimes\phi_i^*=
\tau_1\bold I_1\tag 3.16a\\
&<\phi_i,\phi_j>=\alpha \bold I_{k_1}\tag 3.16b
\endalign$$

Similary, after a complex gauge transformation on $\Cal E_2$, we can pick
$\dbar_2$\ and a basis
$\{\rho_{1},\dots,\rho_{k_2}\}$\ for $V_2$\ such that

$$\align
&i\Lambda F_2+\Sigma_{i={1}}^{k_2}\rho_i\otimes\rho_i^*
=\tau_2\bold I_2\tag 3.17a\\
&<\rho_i,\rho_j>=\alpha \bold I_{k_2}\tag 3.17b
\endalign$$

The gauge transformations on $\Cal E_1$\ and $\Cal E_2$\ combine to produce a
gauge transformation taking $\Cal E$\ to an isomorphic holomorphic bundle,
which we again relable as $\Cal E$.  The $\dbar$-operator corresponding to the
holomorphic
structure on $\Cal E$\ then has a block decomposition as
$$\dbare=\pmatrix
\dbar_1& \beta \\
0 & \dbar_2
\endpmatrix \ ,$$
where $\beta\in \Omega^{0,1}(Hom(\Cal E_2,\Cal E_1))$\ is a
representative of the extension class.  We obtain a properly normalized basis
for $V$\ as follows.  We fix a real number $0<\lambda<1$, and pick smooth
sections $\sigma_i\in\Omega^0(X,E_1)$\ such that
$\tilde\phi_i =\sigma_i+\lambda\rho_i$\ is a lift of $\lambda\rho_i$\ to $\Cal
E$, i.e. $\tilde\phi_i\in H^0(X,\Cal E)$.  Furthermore, the $\sigma_i$\ can be
chosen such that
$$\align
&<\phi_i,\sigma_j>=0\tag 3.18a\\
&<\sigma_i,\sigma_j>=
\alpha(1-\lambda^2)\bold I_{k_1}\tag 3.18b
\endalign$$

Then  $\{\phi_1,\dots,\phi_{k_1},
\tilde\phi_{1},\dots, \tilde\phi_{k_2}\}$\ is an orthogonal basis the basis
for
$V$, with all vectors of length $\sqrt{\alpha}$.  Using $\dbare$\ and this
basis for $V$, we
thus get the block decomposition
$$\align
&i\Lambda F +\Sigma_{i=1}^{k}\phi_i\otimes\phi_i^* \\
&=\pmatrix
\tau_1\bold I_1 -i\Lambda\beta\wedge\beta^* +
\Sigma_{i=1}^{k_2}\sigma_i\otimes\sigma_i^* & B\\
B^* &
\tau_2\bold I_2 -i\Lambda\beta^*\wedge\beta
+(\lambda^2-1)\Sigma_{i=1}^{k_2}\rho_i\otimes\rho_i^*
\endpmatrix \ ,
\endalign$$
where
$$B=\dbar_{1,2}^*\beta + \lambda
\Sigma_{i,j} \sigma_i\otimes\tilde\rho_j^*\ .$$
After a gauge transformation of the form
$\pmatrix
1 & u\\
0 & 1
\endpmatrix$,  the term $B$ changes to
$$B^{(u)}=B+\Delta(u)+
u\lambda^2\Sigma_{i=k_1+1}^{k}\rho_i\otimes\rho_i^*\ ,$$
where $\Delta$ is the ``$\dbar_{1,2}$-Laplacian".
In fact, the endomorphism $u\in\Omega^0(Hom(E_2,E_1)$\ can
be chosen so that $B^{(u)}=0$. This can be seen as follows:
The operator $\Delta
+\lambda^2\Sigma_{i=k_1+1}^{k}\rho_i\otimes\rho_i^*$\ is elliptic,
and its kernel can be identified with those morphisms
$u:(\Cal E_2,V_2)\-->(\Cal E_2,V_2)$\ such that $u(V_2)=0$.
In view of our assumptions on $(\Cal E_2,V_2)$\ and $(\Cal
E_2,V_2)$, it follows by Lemma 2.8 that this kernel is
trivial. The equation $B^{(u)}=0$\ can thus be solved for
$u$.

Suppose then that $B=0$.  The holomorphic structure on
$\Cal E$\ is thus determined by $(\dbar_1,\dbar_2,\beta)$.  Now consider the
1-parameter family
$\{(g_s,A_s)\}$\ in $\Gc(E)\times \GL(k,\Bbb C)$, where
$$g_s=\pmatrix 1 && 0\\
0 && s^{-1}\endpmatrix\ , A_s=\pmatrix 1 && 0\\
0 && \gamma(s) \endpmatrix, \text{and}\
\gamma(s)= \frac{s}{\sqrt{s^2(1-\lambda^2)+\lambda^2}}\ .$$
This generates a 1-parameter family, say $(\Cal E_s,V_s)$,
in the $\Gc(E)\times \GL(k,\Bbb C)$\ orbit through
$(\Cal E, V)$.  The holomorphic structure of $ \Cal E_s$\ is given by
$(\dbar_1,\dbar_2,s\beta)$\
and a basis for $V_s$\ is given by
$\{\phi_1,\dots,\phi_{k_1},\tilde\phi^{(s)}_1,
\dots,\tilde\phi^{(s)}_{k_2}\}$, where
$$\tilde\phi^{(s)}_{k_2}= \gamma(s)\sigma_i+\frac{\lambda}{s}\gamma(s)\rho_i\
.$$

\noindent Observe that $\{\phi_1,\dots,\phi_{k_1},\tilde\phi^{(s)}_1,
\dots,\tilde\phi^{(s)}_{k_2}\}$\ is an orthogonal basis in which all vectors
have length $\sqrt{\alpha}$. Also, in the limit $s\--> 0$, we have
$\gamma(s)\--> 0$\ and $\frac{\lambda}{s}\gamma(s)\--> 1$.  We thus get

$$i\Psi^{CS}_{\G}(\Cal E_s,V_s)-\tau\bold I =
\pmatrix
(\tau_1-\tau)\bold I_1 + \epsilon_1(s)& 0\\
0 &
(\tau_2-\tau)\bold I_2 -\epsilon_2(s)
\endpmatrix \ ,
$$
where
$$\epsilon_1(s)= -s^2i\Lambda\beta\wedge\beta^* +
\gamma(s)^2\Sigma_{i=1}^{k_2}\sigma_i\otimes\sigma_i^*\ ,$$
$$\epsilon_2(s)= s^2i\Lambda\beta^*\wedge\beta
+(1-\frac{\lambda^2}{s^2}\gamma(s)^2)\Sigma_{i=1}^{k_2}\rho_i\otimes\rho_i^*\
.$$

Clearly $\epsilon_1(s)\-->0,\ \epsilon_2(s)\-->0$\ as $s\-->0$. Also, since
$(\Cal E,V)$\ is $\alpha$-stable, we have $\tau_1<\tau<\tau_2$. The proof of
Lemma 3 in [Do] is now easily adapted to show
$$J(\Cal E_s,V_s)=J_1-2s^2|\beta|^2-
k(1-\lambda^2)\gamma(s)^2\alpha\ ,$$
where
$$J_1=R_1(\mu_{\alpha}(\Cal E ,V )-\mu_{\alpha}(\Cal E_1,
V_1))
+R_2(\mu_{\alpha}(\Cal E_2 ,V_2) -\mu_{\alpha}(\Cal E, V))
\ .$$
Thus for small enough $s$, $J(\Cal E_s,V_s)<J_1$, as required.

If $(\Cal E_1,V_1)$\ and $(\Cal E_2,V_2)$\ are not
$\alpha$-stable, then by Lemmas 2.9 and 2.10 they have filtrations by
$\alpha$-stable coherent systems. Exactly as in [Do], the
above arguement can be modified to take this into account. This completes the
proof of Lemma 3.9.

\qed\enddemo

We now apply Lemmas 3.8 and 3.9 to the exact sequences in (3.13).  Reasoning
as
in [Do], we get
$$\align
J(\Cal E,V)&<R_K(\mu_{\alpha}(\Cal E ,V )-
\mu_{\alpha}(\Cal K, V_K))
+R_I(\mu_{\alpha}(\Cal I ,V/V_K)
 -\mu_{\alpha}(\Cal E, V))\\
&\le R_Q(\mu_{\alpha}(\Cal E_{\infty},V_{\infty})-
\mu_{\alpha}(\Cal Q, V_Q))
+R_{\tilde I}(\mu_{\alpha}(\tilde{\Cal I} ,\tilde
g_{\infty}(V))
 -\mu_{\alpha}(\Cal E_{\infty},V_{\infty}))\\
&\le J(\Cal E_{\infty},V_{\infty})\ .
\endalign$$
This is impossible, since $(\Cal E_{\infty},V_{\infty})$\ is minimizing for
$J|_{\Cal O(\E,V)}$.  We have thus proven

\proclaim{Proposition 3.10} If $(\E, V)$\ is an $\alpha$-stable coherent
system, then $(\E_{\infty}, V_{\infty})$\ is on the same $\Gc$-orbit as
$(\E,V)$.
\endproclaim

%%%%%%%%%%%%%%%%%%%%%%%%%%%%
\subheading{\hfil Step 3: At its minimum $J=0$ }
%%%%%%%%%%%%%%%%%%%%%%%%%%%%

  Finally, we must show that by minimizing $J|_{\Cal
O(\E,V)}$\ we obtain a unique solution
to the orthonormal vortex equations. Here it is convenient
to use the smooth function $||\Psi^{CS}_{\G}(\Cal
E,V)+i\tau\I||^2_{L^2}$, rather than the function
$J(\Cal E,V)=N(\Psi^{CS}_{\G}(\Cal E,V)+i\tau\I)$.  Because
of the equivalence of the usual $L^2$\ norm and the norm
used to define $J$, we have
$$C_1 ||\Psi^{CS}_{\G}(\Cal E,V)+i\tau\I||^2_{L^2}
\le   J (\Cal E,V)
\le C_2 ||\Psi^{CS}_{\G}(\Cal E,V)+i\tau\I ||^2_{L^2} \
,$$
for some fixed constants $C_1,C_2$.  It will thus suffice to
show that
$||\Psi^{CS}_{\G}(\Cal E,V)+i\tau\I ||^2_{L^2} $\ has a
unique smooth minimum on the orbit $\Cal O(\E,V)$, and that
its minimum value is zero.

\proclaim{Lemma 3.11} Let $ \G^*=\G/S^1$, where $S^1$\ is
identified with
the subgroup of $\G$\ consisting of the constant multiples
of the
identity.  Then $\G^*$\ acts symplectically on $\Cal
H^{CS}$, and this
action
extends to a holomorphic action of $(\Gc)^*=\Gc/\Bbb C^*$.
The moment
map
for the action of $ \G^*$\ is $\Psi^{CS}_{\G}(\Cal
E,V)+i\tau\I$.
\endproclaim

\demo{Proof} Everything except the last sentence is
immediate. The claim
concerning the moment map follows from the fact that
$\int_X \Tr(\Psi^{CS}_{\G}(\Cal E,V))=-i\tau$. Thus
$\Psi^{CS}_{\G}(\Cal
E,V)+i\tau\I$\ is the projection of $\Psi^{CS}_{\G}(\Cal
E,V)$\ onto the
(dual
of)the Lie algebra of $\G^*$, as required.
\qed\enddemo

It thus follows by general properties of moment maps in such
situations (cf. [K]), that if $\G^*$\ acts freely (or with finite
stabilizer subgroups) on all points of the orbit $\Cal
O(\E,V)$, then $\Psi^{CS}_{\G}(\Cal E,V)+i\tau\I=0$\ at a
minimum for $||\Psi^{CS}_{\G}(\Cal E,V)+i\tau\I||^2_{L^2}$
on the orbit.  Furthermore, there is at most one such
minimum on each orbit.

\proclaim{Lemma 3.12} If $p=(\Cal E,V)$\ is an $\alpha$-stable point in
$\Cal H^{CS}$,
then the
isotropy subgroup of $\G^*$\ at $p$\ is at most finite.
\endproclaim

\demo{Proof} This follows, in the usual way,  from  the fact
that
($\alpha$-)stable
objects are simple. Thus the only automorphisms of $(\Cal
E,V)$\ are the
constant
multiples of the identity.
\qed\enddemo

It thus follows  from the above that if $(\Cal E,V)$\ is
an $\alpha$-stable coherent system with $\rank(\Cal E)\ge 1$, then on $\Cal
O(\Cal E,V)$\ there is a
unique solution to the equation
$$\Psi^{CS}_{\G}(\Cal E',V')+i\tau\I=0\ .$$
The smoothness of the solution follows from standard
elliptic regularity
arguements. This concludes the proof of Theorem 3.2, and thus establishes

\proclaim{Corollary 3.13 (Theorem B)} If $(\Cal E,V)$\ is an $\alpha$-stable
coherent
system, then
there is a unique smooth metric on $\Cal E$, say $H$,  and
frame
$\{\phi_1,\dots,\phi_k\}$\ for $V$\ such that

$$\align i\Lambda F_{\dbare,H}+
\Sigma\phi_i\otimes\phi_i^*&=\tau\bold I\\
 <\phi_i,\phi_j>&=\alpha\bold I_k\ .
\endalign $$
\endproclaim

\noindent\bf {Remark.}\rm\ In the special case where $E$\ is a line bundle,
$(\Cal E,V)$\ is automatically $\alpha$-stable.  Furthermore, the orthonormal
vortex equations reduce to a set of equations which can be viewed as a
Kazdan-Warner type of equation with constraints. The existence theorem which
follows from the above result is, as far as we know, new even in this context.

%%%%%%%%%%%%%%%%%%%%%%%%%%%%%%%%%%%%%%%%%%%%%%%

\end